
\input amstex.tex
\documentstyle{amsppt}
\magnification1200
\hsize=12.5cm
\vsize=18cm
\hoffset=1cm
\voffset=2cm

\def\DJ{\leavevmode\setbox0=\hbox{D}\kern0pt\rlap
{\kern.04em\raise.188\ht0\hbox{-}}D}
\def\dj{\leavevmode
 \setbox0=\hbox{d}\kern0pt\rlap{\kern.215em\raise.46\ht0\hbox{-}}d}

\def\txt#1{{\textstyle{#1}}}
\baselineskip=13pt
\def\hf{{\textstyle{1\over2}}}
\def\a{\alpha}
\def\d{{\,\roman d}}
\def\e{\varepsilon}
\def\f{\varphi}
\def\G{\Gamma}
\def\k{\kappa}
\def\s{\sigma}
\def\t{\theta}
\def\={\;=\;}

\def\zt{\zeta(\hf+it)}

\def\D{\Delta}
\def\no{\noindent}
\def\R{\Re{\roman e}\,} 
\def\z{\zeta}

 \def\t{\theta}
\def\hf{{\textstyle{1\over2}}}
\def\txt#1{{\textstyle{#1}}}
\def\f{\varphi}

\font\tenmsb=msbm10
\font\sevenmsb=msbm7
\font\fivemsb=msbm5
\newfam\msbfam
\textfont\msbfam=\tenmsb
\scriptfont\msbfam=\sevenmsb
\scriptscriptfont\msbfam=\fivemsb
\def\Bbb#1{{\fam\msbfam #1}}

\def \CC {\Bbb C}

\def \ZZ {\Bbb Z}

\font\ff=cmr8
\def\txt#1{{\textstyle{#1}}}
\baselineskip=13pt

\font\teneufm=eufm10
\font\seveneufm=eufm7
\font\fiveeufm=eufm5
\newfam\eufmfam
\textfont\eufmfam=\teneufm
\scriptfont\eufmfam=\seveneufm
\scriptscriptfont\eufmfam=\fiveeufm
\def\mathfrak#1{{\fam\eufmfam\relax#1}}

\font\tenmsb=msbm10
\font\sevenmsb=msbm7
\font\fivemsb=msbm5
\newfam\msbfam
     \textfont\msbfam=\tenmsb
      \scriptfont\msbfam=\sevenmsb
      \scriptscriptfont\msbfam=\fivemsb
\def\Bbb#1{{\fam\msbfam #1}}

\def \CC {\Bbb C}

\def \ZZ {\Bbb Z}

  \def\rightheadline{{\hfil{\ff
Subconvexity for the Riemann zeta-function and the divisor
problem}\hfil\tenrm\folio}}

  \def\leftheadline{{\tenrm\folio\hfil{\ff
   M.N. Huxley and A. Ivi\'c }\hfil}}
  \def\emptyheadline{\hfil}
  \headline{\ifnum\pageno=1 \emptyheadline\else
  \ifodd\pageno \rightheadline \else \leftheadline\fi\fi}

\font\ff=cmr8
\font\teneufm=eufm10
\font\seveneufm=eufm7
\font\fiveeufm=eufm5
\newfam\eufmfam
\textfont\eufmfam=\teneufm
\scriptfont\eufmfam=\seveneufm
\scriptscriptfont\eufmfam=\fiveeufm
\def\mathfrak#1{{\fam\eufmfam\relax#1}}

\font\tenmsb=msbm10
\font\sevenmsb=msbm7
\font\fivemsb=msbm5
\newfam\msbfam
\textfont\msbfam=\tenmsb
\scriptfont\msbfam=\sevenmsb
\scriptscriptfont\msbfam=\fivemsb
\def\Bbb#1{{\fam\msbfam #1}}

\def \CC {\Bbb C}

\def \ZZ {\Bbb Z}

\def\D{\Delta}
\def\a{\alpha}
 \def\e{\varepsilon}
\def\no{\noindent} \def\d{\,{\roman d}}
\topmatter
\title
Subconvexity for the Riemann zeta-function and the divisor problem
\endtitle
\author
Martin N. Huxley and Aleksandar Ivi\'c
\endauthor

\dedicatory
Dedicated to Professor Matti Jutila on the occasion of his retirement
\enddedicatory
\address
School of Mathematics, University of Cardiff, 23 Senghenydd Road,
Cardiff CF2 4AG, Great Britain
\medskip
Katedra Matematike RGF-a, Universitet u Beogradu,  \DJ u\v sina 7,
11000 Beograd, Serbia
\medskip
\endaddress
\keywords The Riemann zeta-function, subconvexity, the divisor
problem, mean square of $|\zt|$, exponent pairs, Bombieri--Iwaniec method
\endkeywords
\subjclass 11 M 06, 11 N 37
\endsubjclass
\email {\tt huxley\@cardiff.ac.uk,  ivic\@rgf.bg.ac.yu}
\endemail
\abstract
A simple proof of the classical subconvexity bound
$\zt \ll_\e t^{1/6+\e}$ for the Riemann zeta-function
is given, and estimation by more refined techniques is discussed.
The connections between the
Dirichlet divisor problem and the mean square of $|\zt|$ are
analysed.
\endabstract
\medskip
\sevenbf
Bulletin CXXXIV de l'Acad\'emie Serbe des Sciences et des
Arts - 2007, Classe des Sciences math\'ematiques et naturelles,
Sciences math\'ematiques No. 32, pp. 13-32.
\medskip
\endtopmatter
\document

\head
1. Convexity for the Riemann zeta-function
\endhead

 Let as usual
$$
\z(s) \= \sum_{n=1}^\infty n^{-s} = \prod_{p \;\roman {prime}}(1-p^{-s})^{-1}
\qquad(\R s > 1)\leqno(1.1)
$$
denote the Riemann zeta-function (see the monographs [22], [23], and [41]
for an extensive account), defined by analytic continuation for values of
$s = \s+it$ for which $\s\le1$. The zeta-function occupies a fundamental
place in analytic number theory, its connection with prime numbers
being evident from~(1.1). The problem of the estimation of $\zt$ and
the evaluation of its mean square is a central one in
zeta-function theory.
The functional equation for $\z(s)$ is
$$
  \zeta(s) = \chi(s)\zeta(1 - s), \;\chi(s) = 2^s\pi^{s-1}\sin\bigl(
  {{\pi s}\over 2}\bigr)\Gamma(1 - s)\quad(\forall s\in\CC), \leqno{(1.2)}
  $$
and it provides a weak form (a stronger form is with logarithms instead of
the ``$\e$" factors) of the so-called ``convexity" bound  in the
so-called ``critical strip" $0\le\s\le1$. Namely the defining series (1.1)
for $\z(s)$ is convergent for $\s = 1+\e$, and by (1.2) it follows that
$$
\z(-\e+it) \ll_\e |\chi(-\e+it)| \ll_\e |t|^{1/2+\e},\leqno(1.3)
$$
since by Stirling's formula for the gamma-function one has, for fixed $\s$,
$$
\chi(s) = (2\pi/t)^{\s+it-1/2}{\roman e}^{i(t+\pi/4)}\Bigl(1+ O(t^{-1})\Bigr)
\qquad(t \ge t_0 > 0).\leqno(1.4)
$$
As usual $\e\,(>0)$ denotes arbitrarily small constants, not necessarily
the same ones at each occurrence, while $a \ll_\e b$ means that the implied
$\ll$--constant depends (only) on $\e$. From $\z(1+\e+it) \ll_\e 1$, (1.3)
and convexity, namely the Phragm\'en--Lindel\"of principle (e.g., see
the Appendix to~[23]), one obtains the desired convexity bound
$$
\z(\s+it) \ll_\e t^{{1\over2}(1-\s)+\e}
   \qquad (0\le \s \le 1,\,t\ge 2).\leqno(1.5)
$$
As mentioned, one can replace ``$\e$" in~(1.5) by a log-power by
using the easily established bounds
$$
\z(1+it) \ll \log t,\qquad \z(it) \ll \sqrt{t}\log t\qquad(t\ge 2).
$$
Naturally, one wishes to obtains sharper bounds than~(1.5),
namely bounds where the exponent of $1-\s$ in~(1.5) is less
than~1/2 for $\s\ge1/2$. The first such exponent
was 1/3 (for $\hf \le \s \le 1$), obtained by Hardy and Littlewood
in~1921 (see Notes to [23, Chapter 7] for a historic discussion of
bounds for~$\zt$), and the best possible bound, up to~``$\e$",
is $\z(\s+it) \ll_\e t^\e$ for $\s \ge 1/2$, which is known as the
Lindel\"of hypothesis. The best known result, that
$\zt \ll_\e t^{32/205+\e},\, 32/205 = 0.15609\ldots\,$
is due to M.N. Huxley~[14].
This is the result of almost a century of continued research,
obtained by varied and refined techniques,
and the current exponent 32/205 is obviously very far from the
exponent zero suggested by the Lindel\"of hypothesis.
This shows the great difficulty of proving the Lindel\"of hypothesis,
and indirectly indicates the enormous difficulty of the Riemann
hypothesis (that all complex zeros of $\z(s)$
have real part 1/2). Namely it is not very difficult
to show (see [36], [23, Chapter 1]) that the Riemann hypothesis,
considered by many to be the greatest open problem
in Mathematics, implies the Lindel\"of hypothesis
(it is not known whether the converse implication is true).
In fact if the Riemann hypothesis is true,
then (see E. C. Titchmarsh's monograph [41]) one has the bound
$$
\zt \;\ll \; \exp\left(A{\log t\over\log\log t}\right)
\qquad(A>0, t \ge 2),
$$
which is stronger than the Lindel\"of hypothesis,
namely the bound $\zt \ll_\e t^\e$.

\medskip
There is a natural, very important aspect of this subject which
concerns ``convexity" bounds for a large class of Dirichlet series.
This is the so-called Selberg class $\Cal S$, consisting of
Dirichlet series $F(s) = \sum_{n=1}^\infty f(n)n^{-s}\;(\s>1)$
which possess several properties analogous to $\z(s)$,
the most important ones being the Euler product over primes
(see (1.1)) and the functional equation (see (1.2)) involving
gamma-factors. For a comprehensive survey of $\Cal S$, the reader is
referred to the paper of J.~Kaczorowski and A.~Perelli [34].
For each function $F\in {\Cal S}$ there exists a convexity bound
analogous to~(1.5). In this general context `subconvexity' means a
bound which improves on the appropriate analogue of~(1.5),
especially when $\s = \hf$, which is the so-called ``critical line"
in the theory of functions in $\Cal S$.

\medskip
The aim of this paper is twofold. Firstly, we shall provide a simple,
self-contained subconvexity bound for $\z(s)$, the prototype of all
functions from $\Cal S$. We shall also give a brief outline of the
Bombieri--Iwaniec method at the end of the paper,
which furnished the hitherto best bounds for~$\zt$.
Secondly, we shall analyse the connection between the mean square
of $|\zt|$ and the classical Dirichlet divisor problem,
namely the estimation of the function
$$
\D(x) = \sum_{n\le x}d(n) - x(\log x + 2\gamma -1),\leqno(1.6)
$$
where $d(n)$ is the number of divisors of~$n$ and
$\gamma = - \G'(1) = 0.57721\ldots\,$
is Euler's constant. This connection is a natural one, since
from~(1.1) it follows that
$$
\z^2(s) \= \sum_{n=1}^\infty d(n)n^{-s}\qquad(\s>1),
$$
so that $d(n)$ is generated by $\z^2(s)$. It is not yet clear whether
an analogue of such a connection exists in general for elements
of~$\Cal S$, or represents an intrinsic property of a subclass
of~$\Cal S$ containing $\z(s)$.

\head
2. Subconvexity for the Riemann zeta-function
\endhead
We start from [22, Theorem 1.2] or [23, Lemma 7.1], which bounds
$\z(\hf+iT)$ by its mean square over a short interval:
$$\eqalign{
|\z(\hf+iT)|^2 &\ll
 \log T\left(1+\int_{T-\log^2T}^{T+\log^2T}|\zt|^2\d t\right)\cr
&\ll \log T\int_{T-2G}^{T+2G}f(t)|\zt|^2\d t.\cr}\leqno(2.1)
$$
Here $T^\e\le G\le T^{1-\e}$, and $f(t)$ is a smooth, non-negative
function supported in $[T-2G,T+2G]$ such that $f(t) = 1$ for
$T-G\le t \le T+G$.
The second bound in~(2.1) is trivial, while the first one rests on
contour integration of $\z^2(s+z)\G(z)$ over $z$ and the
use of the functional equation (1.2) for~$\z(s)$,
and thus may be considered to be elementary. Although (2.1) appears
to be wasteful, it turns out that it is an effective starting point
for the estimation of~$\zt$.

By the approximate functional equation for $\z(s)$
(see e.g., Theorem~4.4 of~[23]) the second integral
in~(2.1) is majorized by $O(\log T)$ integrals of the type
$$\eqalign{
I :&= \int_{T-2G}^{T+2G}f(t)
  \Bigl|\sum_{N<n\le N_1}n^{-{1\over2}-it}\Bigr|^2\d t\cr
&= O(G)  + \sum_{N<m\ne n\le N_1}(mn)^{-{1\over2}}
 \int_{T-2G}^{T+2G}f(t)(m/n)^{it}\d t,\cr}\leqno(2.2)
$$
where $1\ll N < N_1 \le 2N \ll T^{1/2}$.
We may assume $G \ll N$, for otherwise the contribution of~$I$ is
$\ll G$ by the well-known mean value theorem for Dirichlet
polynomials (see e.g., Theorem 5.2 of~[23]).
By symmetry it may be assumed that $m>n$, thus $m = n+r,\,r\ge1$.
Note that
$$
\int_{T-2G}^{T+2G}f(t) (m/n)^{it}\d t = {i\over\log{m\over n}}
\int_{T-2G}^{T+2G}f'(t) (m/n)^{it}\d t,
$$
so that integrating
the last integral in~(2.2) sufficiently many times by parts and
using $f^{(j)}(t)\ll_j G^{-j},\,1/\log(m/n) \ll N/r$,
it follows that the contribution of $r > N^{1+\e}/G\;$
is~$\;\ll 1$. Therefore we have
$$
I \;\ll_\e\;
  G + \Biggl|\sum_{1\le r\le N^{1+\e}/G}\int_{T-2G}^{T+2G}f(t)
  \sum_{N<n\le N_1}n^{-{1\over2}}(n+r)^{-{1\over2}}
  {(1+r/n)}^{it}\d t\Biggr|.\leqno(2.3)
$$
The sum over~$n$ in~(2.3) is written as
$$
\sum_{N-G<n\le N_1+G}\f(n)n^{-1/2}(n+r)^{-1/2}{(1+r/n)}^{it}
  + O(G/N),\leqno(2.4)
$$
where $\f(t)\,(\ge0)$ is a smooth function supported
in $[N-G,N_1+G]$ such that $\f(t) = 1$ for $N\le t \le N_1$,
and thus $\f^{(j)}(t) \ll_j G^{-j}$. The sum in~(2.4) is treated
by the Poisson summation formula ($F(x)$ is real valued,
smooth and compactly supported in $[0,\infty)$)
$$
\sum_{n=1}^\infty F(n)
  = \int_0^\infty F(x)\d x + 2\sum_{k=1}^\infty
  \int_0^\infty F(x)\cos(2\pi kx)\d x,
$$
applied with $F(x)$ equal to the real (respectively imaginary) part
of $\f(x)x^{-1/2}(x+r)^{-1/2}(1+r/x)^{it}$. Performing again
sufficiently many integrations by parts it is seen that the
contribution of $k > T^{1+\e}N^{-2}r\;$ will be~$\ll 1$.
We are left with the sums
$$
\sum_{k\le T^{1+\e}N^{-2}r}\int_{N-G}^{N_1+G}
  \f(x)x^{-1/2}(x+r)^{-1/2}(1+r/x)^{it}
  \exp(\pm2i\pi kx)\d x.
$$
The exponential factor here is of the form
$$
\exp(ig(x)),\quad g(x):= t\log\Bigl(1+{r\over x}\Bigr) \pm 2k\pi x,
\quad g''(x) = t{2rx+r^2\over(x^2+rx)^2} \ll {rT\over N^3}.
$$
Thus applying the second derivative test (Lemma~2.2 of~[23]) to
the integral we deduce that it is $\ll N^{-1}(N^3/(rT))^{1/2}$.
By the first derivative test ([23] Lemma~2.1) the total
contribution of $\int_0^\infty F(x)\d x$ is clearly $\ll_\e GT^\e$.
Hence from the preceding estimates we have
$$\eqalign{
I &\ll_\e T^\e\left(G + G\sum_{r\le N^{1+\e}/G}N^{1/2}(rT)^{-1/2}
\sum_{k\le T^{1+\e}N^{-2}r}1\right) \cr
&\ll_\e T^\e(G + T^{1/2}G^{-1/2}) \ll_\e T^\e G
\cr}
$$
for $G \ge T^{1/3}$. Hence, with $G=T^{1/3}$, (2.1) yields the bound
$\zt\ll_\e t^{1/6+\e}$,
which is the desired subconvexity estimate.
As is to be expected, more refined estimates lead to
exponents $< 1/6$, which will be shown in Section~4.

\bigskip
There are other ways to obtain the desired subconvexity estimate.
For example, one can use [23, Theorem 1.8]. For the estimation of
the zeta-sum
$$
\sum_{N<n\le N_1}n^{-1/2-it}
$$
one then uses the exponent pair $({1\over6},\,{2\over3})$
([23] Chapter 2 or M.N. Huxley [8]), which immediately leads to
the classical estimate of Hardy--Littlewood
(see e.g., E.C. Titchmarsh's monograph [41, Chapter 5])
$$
\zt \;\ll\; t^{1/6}\log t.\leqno(2.5)
$$
The fact is that one can show in an elementary way that
$({1\over6},\,{2\over3})$ is an exponent pair,
by using [23, Lemma 2.5] (Weyl's inequality) and then Poisson
summation (in the form of, say, [23, Lemma 7]). This is
probably the quickest way to attain the subconvexity bound~(2.5),
but perhaps not self-contained as the previous approach of ours.
For two related approaches the reader is referred to Chapter~5
of Titchmarsh's book~[41].

\bigskip
One can also use the approach of D.R. Heath-Brown [6],
by using P. Gallagher's inequality (e.g., see H.L. Montgomery
[38, Lemma 1.10])
$$
\int_{-G}^G\Bigl|\sum_{n=1}^\infty a_nn^{-it}\Bigr|^2\d t
\;\ll\; G^2\int_0^\infty
\Biggl|\sum_y^{y{\roman e}^{1/G}}a_n\Biggr|^2{\d y\over y}
\leqno(2.6)
$$
applied to $a_n = n^{-iT}$ if $N < n \le N_1 \le 2N$ and $a_n=0$
otherwise. This easily leads  to
$$
\int_{T-G}^{T+G}|\zt|^2 \d t \ll G\log T\qquad(T^{1/3}
  \ll G \ll T),
$$
giving again (2.5) in view of~(2.1). This is quick indeed,
but requires the knowledge of~(2.6).

\head
3. The Dirichlet divisor problem
\endhead

\bigskip
It is interesting that one can also establish the
analogy of the mean square of $|\zt|$ with the classical
Dirichlet divisor problem, namely the estimation of the
error-term  function~$\D(x)$, defined by~(1.5).
By using the elementary formula (see [23, eq. (14.41)])
$$
\sum_{n\le z}n^{-1} \;=\;
\log z + \gamma - \psi(z)z^{-1} + O(z^{-2})\quad(z\ge2),
\quad \psi(x) = x - [x] - \hf,
$$
and writing
$$
\sum_{n\le x}d(n) = \sum_{mn\le x}1
= 2\sum_{n\le \sqrt{x}}[x/n] - [\sqrt{x}\,]^2,
$$
one easily arrives at
$$
\D(x) \;=\; -2\sum_{n\le \sqrt{x}}\psi(x/n) + O(1).\leqno(3.1)
$$
But, for $T^\e \le G \le T^{1/2}$ and suitable $C>0$, we
trivially have
$$
C\D(T) - {1\over G}\int_{T-2G}^{T+2G}f(t)\D(t)\d t =
{1\over G}\int_{T-2G}^{T+2G}f(t)(\D(T)-\D(t))\d t \ll_\e GT^\e,
$$
since $d(n) \ll_\e n^\e$, where $f(t)$ is as in~(2.1).
Using the Fourier expansion
$$
\psi(x) = -{1\over\pi}\sum_{m=1}^\infty {\sin(2\pi mx)\over m}
\qquad(x\not\in \ZZ)
$$
and the fact that the above series is boundedly convergent,
 it follows from~(3.1) that
$$
\D(T) = {2\over \pi CG}\sum_{n\le \sqrt{T}}\sum_{m=1}^\infty
\int_{T-2G}^{T+2G}f(t){\sin(2\pi mt/n)\over m}\d t + O_\e(GT^\e).
$$
If we perform sufficiently many integrations by parts, using the
fact that $f^{(\ell)}(t) \ll_\ell G^{-\ell}$
for $\ell = 0,1,2,\ldots\,$, it transpires that the contribution
of~$m$ satisfying $m > nT^\e G^{-1}$ is negligible, that is,
it is $O(1)$. Hence we are left with the estimate
$$\eqalign{
\D(T) &\ll {\log^2T\over G}\sup_{M\le N^{1+\e}/G,N\ll\sqrt{T}}
\int\limits_{T-2G}^{T+2G}f(t)\Bigl|{1\over M}\sum_{M<m\le M'}
\sum_{N<n\le N'}\exp\Bigl({2\pi i mt\over n}\Bigr)\Bigr|\d t \cr
&+ GT^\e,\cr}\leqno(3.2)
$$
where $M < M'\le 2M,\, N< N'\le 2N$. The bound in~(3.2) is the
analogue of~(2.3), if we notice that in~(2.3) we have
$$
\left(1 + {r\over n}\right)^{it}
= \exp\left\{it\log\left(1 + {r\over n}\right)\right\}
= \exp\left\{it\left({r\over n}
   - {1\over2} \left({r\over n}\right)^2
   + \ldots\, \right)\right\},\leqno(3.3)
$$
so that the term $itr/n$ dominates in the exponential in~(3.3).
Therefore we easily obtain from~(3.2) the bound
$\D(x) \ll_\e x^{1/3+\e}$.
A short proof of the slightly sharper classical bound
$\D(x) \ll x^{1/3}\log x$, by the use of the Vorono{\"\i} formula
(see (5.6)), is given by the author in~[24], but the Vorono{\"\i}
formula may be considered as a non-elementary tool, while the
preceding discussion was completely elementary.

\head
4. Sharper bounds
\endhead
\medskip
It is worth remarking that both in the case of~$\D(x)$ and
in the case of the mean square of~$|\zt|$ one can
obtain notably better bounds. For the latter, as usual, we define
$$
E(T) \;=\;\int_0^T|\zt|^2\d t
 - T\left(\log\bigl({T\over2\pi}\bigr)   + 2\gamma - 1 \right)
$$
and note that we have the elementary inequalities
(see the author's monograph [22, Lemma~4.1])
$$
E(T) \le E(T+x) + Cx\log T,\quad E(T) \ge E(T-x)
  - Cx\log T \quad(0\le x\le T). \leqno(4.1)
$$
We replace $T+x$ by $t$, multiply the first inequality in~(4.1)
by $f_0(t)$ and integrate, where $f_0(t)\;(\ge0)$ is a smooth
function supported in $[T,\,T+3G]$ such that $f_0(t) = 1$ when
$T+G\le t \le T+2G$.  It follows that, with suitable $C_1,\,C_2$,
$$
E(T) \le {C_1\over G}\int_{T}^{T+3G}f_0(t)E(t)\d t + C_2G\log T
\quad(C_1 >0,C_2 >0, \,T^\e \le G \le \sqrt{T}\,),\leqno(4.2)
$$
and also an analogous lower bound inequality holds. To deal with
the integral in~(4.2) we invoke the explicit formula of
R.~Balasubramanian~[2], whose work is based on the use of the
classical Riemann--Siegel formula for $\zt$ (see [23, eq. (4.5)]
or [41, Theorem 4.16]). This is
$$
\eqalign{
E(T) &=  2\sum_{n\le K}\sum_{m\le K,m\ne n}
  {\sin(T\log n/m)\over\sqrt{mn}\log n/m} \cr
  &+ 2\sum_{n\le K}\sum_{m\le K,m\ne n}
  {\sin(2\t_1-T\log mn)\over\sqrt{mn}\,(2\t'_1-\log mn)}
  + O(\log^2T),
\cr}\leqno(4.3)
$$
where
$\t_1 = \t_1(T) = \hf T\log (T/(2\pi)) - \hf T - {1\over8}\pi,
K = \sqrt{T/(2\pi)}$. We insert (4.3) (with $T$ replaced by $t$)
in~(4.2), and integrate by parts as before, using
$f_0^{(j)}(t) \ll_j G^{-j}$.
The double sums in~(4.3) will yield exponential factors of the form
$$
\exp\left(it\log{m\over n}\right),\quad
\exp\left(it\log{tmn\over 2\pi} - it\right)\quad(m\ne n),\leqno(4.4)
$$
and the sums over $m$ and $n$ are split in $\ll \log^2T$ subsums with
$$
M < m \le M' \le 2M,\quad N < n \le N' \le 2N,
\quad M \ll \sqrt{T},\, N\ll \sqrt{T}.
$$
Note that
$$
{\d\over \d t}\left(t\log{m\over n}\right) = \log{m\over n},\;
{\d\over \d t}\left(t\log{tmn\over 2\pi} - t\right)
   = \log{tmn\over 2\pi}.
$$
Thus, by the analysis that led to~(2.3),
the contribution of the second sum in (4.3) is clearly negligible.
For the first exponential factor in~(4.4)
the non-negligible contribution will come from the values of~$M,N$
for which $M \ll N \ll M$.  Further we set $n = m+ r$ where
$r\ne 0$ is an integer, and it follows that $M \le N^{1+\e}G^{-1}$
is the range for which the contribution is non-negligible.
In view of~(3.3) it transpires that an expression analogous to the
right-hand side of~(3.2) is obtained.

\bigskip

Both exponential sums in (2.3) and~(3.2) are double (two-dimensional)
exponential sums, and as such they may be treated by a variety of
techniques and transformations developed for these types of
sums (see e.g., the book of S.W. Graham and G. Kolesnik~[5]).
However, non-trivial results may be obtained if already the sum
over~$n$, say in~(3.2), is estimated as $\ll (xM/N^2)^\k N^\lambda$
by the theory of one-dimensional exponent pairs
(see e.g., [5, Chapter 3] or [23, Chapter 2]),
where $(\k,\lambda)$ is an exponent pair.
Then the right-hand side of~(3.2) is, by trivial estimation,
$$
\ll_\e
T^\e(T^{(\k+\lambda)/2}G^{-\k} + G)
   \ll_\e T^{{\k+\lambda\over2+2\k}+\e}
$$
with $G = T^{(\k+\lambda)/(2+2\k)}\,(\le \sqrt{T}\,)$, since
$0 \le \k \le \hf \le \lambda \le 1$ has to hold for any exponent
pair~$(\k,\lambda)$. Therefore from the preceding discussion we
obtain a proof of the following

\bigskip
THEOREM. {\it If $(\k,\lambda)$ is an exponent pair, then}
$$
E(T) \ll_\e T^{{\k+\lambda\over2+2\k}+\e},\qquad
\D(T) \ll_\e T^{{\k+\lambda\over2+2\k}+\e}.\leqno(4.5)
$$

\bigskip
From (2.1) and (4.5) we readily obtain then the following

\bigskip
{\bf Corollary}.
$$
\zt \;\ll_\e\; |t|^{{\k+\lambda\over4+4\k}+\e}.\leqno(4.6)
$$

\medskip
With the standard (elementary) exponent pair $(\hf,\,\hf)$ we
obtain again the exponent $1/3+\e$ in~(4.5).
The exponent is less than~1/3 if
$$
3\lambda + \k < 2.\leqno(4.7)
$$
If we take, e.g., the exponent pair $({11\over30},\,{16\over30})$,
then we obtain the exponent $27/82+\e$ ($< 1/3$) in~(4.5).
Note that in~(4.6) any non-trivial exponent pair (i.e., any pair
$(\k,\,\lambda)$ with $\lambda < 1$) beats convexity,
namely produces the exponent in~(4.6) that is strictly less
than~1/4, while any exponent pair satisfying (4.7) improves (2.5).
For some exponent pairs, obtained by the Bombieri--Iwaniec method,
the reader is referred to M.N. Huxley~[12], whilst the method itself
will be analysed in Section~6.
The results in the book~[12] give the exponent $45/137+\e$ in~(4.5),
and there will be a small improvement from the latest work~[14],
and a better improvement from Sargos's programme, but the sum
in~(4.5) is precisely the one that the two-variable form of the
Bombieri--Iwaniec method is designed to treat.

\head
5. Connections between $E(T)$ and $\D(x)$
\endhead

We conclude our discussion by analyzing more closely
the connection between $E(T)$ and $\D(x)$. This has become well-known
after the pioneering work of F.V. Atkinson~[1], who established an
explicit formula for $E(T)$ (different from~(4.3)).
Atkinson's result is the following: let $0 < A < A'$ be any two
fixed constants such that $AT < N < A'T$, and let
$N' = N'(T) = T/(2\pi) + N/2 - (N^2/4+ NT/(2\pi))^{1/2}$. Then
$$
E(T) = \Sigma_1(T) + \Sigma_2(T) + O(\log^2T),\leqno(5.1)
$$
where
$$
\Sigma_1(T) = 2^{1/2}(T/(2\pi))^{1/4}\sum_{n\le N}(-1)^nd(n)n^{-3/4}
e(T,n)\cos(f(T,n)),\leqno(5.2)
$$
$$
\Sigma_2(T) = -2\sum_{n\le N'}d(n)n^{-1/2}(\log T/(2\pi n))^{-1}
   \cos\left(T\log \Bigl( {T\over2\pi n}\Bigr) - T
   + {\txt{1\over4}}\pi\right),
\leqno(5.3)
$$
with
$$
\eqalign{
\cr
&f(T,n) = 2T{\roman {ar\,sinh}}\,\bigl(\sqrt{\pi n/(2T)}\,\bigr)
  +  \sqrt{2\pi nT + \pi^2n^2} - {\txt{1\over4}}\pi\cr
&=  -\txt{1\over4}\pi + 2\sqrt{2\pi nT} +
  \txt{1\over6}\sqrt{2\pi^3}n^{3/2}T^{-1/2} + a_5n^{5/2}T^{-3/2}
  + a_7n^{7/2}T^{-5/2} + \ldots\,,
\cr}\leqno(5.4)
$$
$$\eqalign{\cr
e(T,n) &= (1+\pi n/(2T))^{-1/4}{\Bigl\{(2T/\pi n)^{1/2}
{\roman {ar\,sinh}}\,\Bigl(\sqrt{\pi n/(2T)}\,\Bigr)\Bigr\}}^{-1}\cr
&= 1 + O(n/T)\qquad(1 \le n < T),
\cr}\leqno(5.5)
$$
and $\,{\roman{ar\,sinh}}\,x = \log(x + \sqrt{1+x^2}\,)$.
For $\D(x)$ we have the explicit, truncated Vorono{\"\i}
formula (see e.g., [23] or [42])
$$
\D(x) = {1\over\pi\sqrt{2}}
x^{1\over4}\sum_{n\le N}d(n)n^{-{3\over4}}\cos(4\pi\sqrt{nx}
- {\txt{1\over4}}\pi) +
O_\e(x^{{1\over2}+\e}N^{-{1\over2}})\quad(2 \le N \ll x).
\leqno(5.6)
$$
A comparison between (5.2) and (5.6) reveals at once the
similarities between $E(T)$ and $\D(x)$.
This becomes even more pronounced if one considers
$$
\D^*(x) \;:=\; -\D(x)  + 2\D(2x) - \hf\D(4x)\leqno(5.7)
$$
instead of $\D(x)$. Then  the arithmetic interpretation
of~$\D^*(x)$ (see T. Meurman [37]) is
$$
\hf\sum_{n\le4x}(-1)^nd(n) \;=\; x(\log x + 2\gamma - 1)
   + \D^*(x).
\leqno(5.8)
$$
One also has a Vorono{\"\i}-type formula
(see e.g., [23, eq. (15.68)]), for $2 \le N \ll x$,
$$
\D^*(x) = {1\over\pi\sqrt{2}}x^{1\over4}
\sum_{n\le N}(-1)^nd(n)n^{-{3\over4}}
\cos(4\pi\sqrt{nx} - {\txt{1\over4}}\pi) +
O_\e(x^{{1\over2}+\e}N^{-{1\over2}}),
\leqno(5.9)
$$
which is completely analogous to~(5.6). M. Jutila,
in his works [28] and [29], investigated both the local and
global behaviour of the difference
$$
E^*(t) \;:=\; E(t) - 2\pi\D^*\bigl({t\over2\pi}\bigr).
$$
He proved the mean square bound
$$
\int_{T-H}^{T+H}(E^*(t))^2\d t \ll_\e  HT^{1/3}\log^3T + T^{1+\e}
\quad(1 \ll H \le T), \leqno(5.10)
$$
which in particular yields
$$
\int_0^T(E^*(t))^2\d t \ll T^{4/3}\log^3T.\leqno(5.11)
$$
One conjectures that $\D(x)\ll_\e x^{1/4+\e}$ and
$E(T) \ll_\e T^{1/4+\e}$  both hold
(the exponent 1/4 is in both cases best possible).
A weaker conjecture, still unproved but supported by the bounds
in~(4.5), is that $\a = \gamma$, where
$$
\a = \inf\{\,a > 0 \;:\; \D(x) \ll x^a\,\},
\quad \gamma = \inf\{\,g > 0 \;:\; E(T) \ll T^g\,\}.\leqno(5.12)
$$
For sharp upper bounds on $\a$ and $\gamma$, obtained by the
so-called Bombieri--Iwaniec method, which are
better than those in~(4.5), the reader is referred to the works
of M.N. Huxley [8], [14]. This method will be briefly described
in Section~6. It may be remarked that the conjectural values
$\a = 1/4$ or $\gamma = 1/4$ cannot be attained even if one
assumes the Riemann hypothesis. This clearly shows the difficulty
of this subject. On the other hand, the strong conjecture that
$(\e\,,\hf+\e)$ is an exponent pair
(this implies the hitherto unproved Lindel\"of hypothesis)
does yield the conjectural values $\a = 1/4$ and $\gamma = 1/4$.

\medskip
For an extensive discussion of $E^*(T)$ see the author's
monographs [22],~[23]. The significance of~(5.11) is that it
shows that $E^*(t)$ is in the mean square sense much smaller than
either $E(t)$ or~$\D^*(t)$. It is expected that $E(t)$
and~$\D^*(t)$ are `close' to one another in order,
but this has never been satisfactorily established, although
M.~Jutila [28] obtained significant results in this direction.
It is also conjectured that
$E^*(T) \ll_\e T^{1/4+\e}, \D^*(x) \ll_\e x^{1/4+\e}$,
which is the analogue of the classical conjecture
$\D(x)\ll_\e x^{1/4+\e}$ in the Dirichlet divisor problem.
Recently Y.-K. Lau and K.-M. Tsang [35]
proved that $\a = \a^*$, where $\a$ is as in~(5.12) and
$$
 \a^* \;=\;
\inf\Bigl\{a^* >0\;:\;\D^*(x) \ll x^{a*}\,\Bigr\}.
$$

\medskip
In the first part of the second author's work~[25] the bound
in~(5.11) was complemented with the new bound
$$
\int_0^T (E^*(t))^4\d t \;\ll_\e\; T^{16/9+\e};\leqno(5.13)
$$
neither (5.11) or (5.13) seem to imply each other.
In the second part of the same work [25], it was proved that
$$
\int_0^T |E^*(t)|^5\d t \;\ll_\e\; T^{2+\e},
$$
and some further results on higher moments of~$|E^*(t)|$ were
obtained as well.
In [26] it was shown that the bound in~(5.11) is of the correct
order of magnitude. The result is the asymptotic formula
$$
\int_0^T (E^*(t))^2\d t
  \;=\; T^{4/3}P_3(\log T) + O_\e(T^{7/6+\e}),\leqno(5.14)
$$
where $P_3(y)$ is a polynomial of degree three in~$y$ with
positive leading coefficient, and all the coefficients may be
evaluated explicitly. The formula (5.14) is the limit of the method.

\head 6.
The Bombieri-Iwaniec method in a nutshell
\endhead
\medskip
The sharpest subconvexity bounds for the Riemann zeta-function are
obtained by a method which has still not found a descriptive name.
The originators prefer `Bombieri-Iwaniec method'; Huxley suggested
`the discrete Hardy-Littlewood method'; and Sargos suggested
`Poisson summation modulo one'. It has been successfully applied
to four types of sum:
$$
S_1   \;=\;   \sum _{n \asymp N} b(n){\roman e}(f(n))\qquad\Bigl({\roman e}(z) \;:=
\;\exp(2\pi iz)\Bigr),
$$
where $b(n)$ are the coefficients of a modular form, $n\asymp N$
means that $N<n\le N'$ with $N< N'\le2N$, and $f(x)$ is a real
function growing faster than $O(N)$,
$$
\eqalign{
S_2   &=   \sum _{n \asymp N} {\roman e}(f(n)), \cr
S_3   &=   \sum _{h \asymp H} \sum _{n \asymp N} {\roman e}(hf(n)), \cr}
$$
and
$$
S_4   =   \sum _{h \asymp H} \sum _{n \asymp N}
{\roman e}\Bigl(f(n+h)-f(n-h)\Bigr).
$$
Putting
$$
f(x)  =   - {t\over 2\pi} \log {x\over N} \leqno (6.1)
$$
makes $S_1$ corresponds to partial sum of the Hecke $L$-series $L(s,f)$,
$S_2$ to those of the Riemann zeta-function, and $S_4$ to the short
interval mean square of the Riemann zeta-function according to~(2.4).
Putting
$$
f(x)   =   \sqrt{R^2-x^2} \leqno (6.2)
$$
makes $S_3$ correspond to the Gauss circle problem of the lattice
points in the circle $m^2+n^2 \leq  R^2$, and putting
$$
f(x) = T/x \leqno (6.3)
$$
makes $S_3$ correspond to the Dirichlet divisor problem according to~(3.2).
The sum
$$
S_5   =   \sum _{h \asymp H} \sum_{k \asymp K} \sum _{n \asymp N}
{\roman e}\Bigl(kf(n+h)-kf(n-h)\Bigr)
$$
could be used to investigate the short interval mean square in the
circle problem, but the traditional Fourier series methods already
give results which are in some sense best possible [11, 39].

The method is bedevilled by technicalities and delicate estimates.
There are two great ideas: the local lattice basis, and the large
sieve. First we divide the sum into short intervals according to
Diophantine approximation to $f'(x)$ in $S_1$ and $S_3$, or
$f''(x)/2$ in $S_2$ and~$S_4$. On a subinterval of~$S_1$ where
$f'(x)$ is close to an integer, the sequence $f(n)$ is slowly
varying modulo one, and we use the Vorono\"{\i}-Wilton summation
formula~[47]. If $f'(x)$ is close to a rational number $a/q$, then the
sequence $f(n)-an/q$ is slowly varying modulo one, and we can use
the Vorono\"{\i}-Wilton formula for the coefficients $b(n){\roman e}(an/q)$
of the modular form twisted by an additive character, or by the
action of the matrix
$$
\left(\matrix a& -\bar q\\
q&\bar a\endmatrix\right),
\leqno (6.4)
$$
where $\bar a$ and $\bar q$ are integers with $a\bar a+q\bar q = 1$,
which correspond to a change of basis in the underlying lattice.
This idea was used by Jutila [30, Theorem~4.7] to obtain the result
$$
\int_0^T\left|L\left(\hf\kappa +it,f\right)\right|^6\d t
  \;\ll_\e\; T^{2+\e},
$$
where $L(s,f)$ is the zeta-function of a holomorphic cusp form of
weight~$\k$ for the full modular group. The method was developed
in~[31] and many subsequent papers.

\medskip
Soon afterwards Bombieri and Iwaniec wrote their great paper [3]
on the size of the Riemann zeta-function on the critical line.
The first idea in~[3]
can be interpreted as the local lattice basis step for the sum
$$
S_2'  \; =\;   \sum _{m \asymp N^2} q(m){\roman e}\bigl(f(\sqrt{m}\,)\bigr),
$$
where $q(m) = 2$ if $m$ is a perfect square, $0$ otherwise, which is
a lacunary sum of the form~$S_1$, with the coefficients $b(m)$
those of a Jacobi theta function.

Iwaniec and C.J.~Mozzochi [27] considered the more complicated sum~$S_3$,
which corresponds to a two-dimensional lattice point problem,
with~$a/q$ a rational approximation to the gradient of
the boundary curve, in the case~(6.3). The sum $S_4$ was treated
by Heath-Brown and Huxley~[7]. The rational number $a/q$ in~(6.4)
is an approximation to~$2f'(x)$, and the local lattice basis step
follows that for~$S_3$.

The local lattice basis step on its own leads to the estimates
$$
\eqalign{
L(\hf\k+it,f)  \;&\ll _\e\;  t^{1/3+\e },\cr
\zeta \left(\hf+it\right)  \;&\ll _\e\;  t^{1/6 + \e }, \cr
\Delta (x)  \;&\ll _\e\;  x^{1/3 + \e}, \cr
E(T) \;&\ll _\e\;  T^{1/3 + \e }. \cr}
$$
Further savings come from estimating the transformed sums after
Vorono\"{\i} or Poisson summation. The transformed sum can be
estimated directly when the denominator $q$ is small (the `major
arc' case).

Bombieri and Iwaniec's second innovation was to use the large
sieve inequality in a general form of their devising. The first
three or four terms of the power series for the function in the
exponent are regarded as a vector inner product between a
coefficient vector which depends only on the short interval, and a
variable vector such as $(h^2,h^{3/2},h,\sqrt{h})$, whose entries
are monomials in the new variables introduced by Vorono\"{\i} or
Poisson summation.

The variable vectors form a thin set, and H\"older's inequality is
used to pass to powers of the transformed sums, so that four or
five or six variable vectors with different integer parameters are
added. Large sieves usually require the vectors to be distinct,
but the Bombieri-Iwaniec form has a sum over a `neighbourhood of
the diagonal', that is, over pairs of vectors which are close
enough together for their inner products with coordinate vectors
to be approximately equal. The First Spacing Problem is to
estimate this sum. Usually the parameters are chosen so that the
order of magnitude of the sum is essentially that of the diagonal
contribution. The First Spacing Problem for the sum~$S_3$ (a sum
of four vectors) was settled by Iwaniec and Mozzochi. For the
sum~$S_2$, Bombieri and Iwaniec [3,4] took a sum of four vectors in
the First Spacing Problem. N.~Watt [43,44] simplified the treatment,
and Watt~[45] and Huxley and Kolesnik~[15]
found partial results for a sum of five vectors; the optimal use
of H\"older's inequality would give a sum of six vectors. The
First Spacing Problem for~$S_4$ is like that for~$S_2$, but with
extra low order terms. Heath-Brown and Huxley~[7] had a partial result,
and Watt~[46] has recently obtained the expected bound.

The Second Spacing Problem is to estimate a sum over pairs of
coefficient vectors which form a neighbourhood of the diagonal.
The coefficients depend on the function~$f(x)$. The problem is
essentially the same for the sums $S_2$, $S_3$, and $S_4$. Huxley
and Watt~[18] and independently Kolesnik (unpublished) generalised
Bombieri and Iwaniec's bound from the special case~(6.1) in~$S_2$ to a
general function~$f(x)$. Huxley~[8], and independently Li
Hongquan (unpublished) generalised the treatment of~$S_3$ from
Iwaniec and Mozzochi's special case~(6.3) to a general
function~$f(x)$. The large sieve is applied to the transformed sum
on the minor arcs, and the conditions for a pair of minor arc
coefficient vectors labelled $a/q$ and $a'/q'$ to lie in a
neighbourhood of the diagonal are expressed in terms of the
lattice base change matrix
$$
M = \left(\matrix a & -\bar q\\q &\bar a\endmatrix\right)
    \left(\matrix a'  &-\bar q' \\
              q'  &\bar a' \endmatrix\right)^{-1}.
\leqno (6.5)
$$
The basic counting idea is that $M$ acts only for short ranges of
the rationals $a/q$ and~$a'/q'$; this idea uses only two of the
four entries of the coefficient vector. Huxley [10,14] has small
improvements in which some use is made of the other two entries of
the coefficient vector. For a fixed base change matrix $M$,
the other two conditions can be interpreted as saying that an integer
point in some dual plane lies close to a curve determined by the
matrix~$M$. The parameters are chosen so that the order of
magnitude of the sum is essentially that of the diagonal
contribution.

Jutila~[32,33], treating the sum~$S_1$, introduced a variant where
the minor arcs overlap, and so there are more denominators $q$ in
the Second Spacing Problem.

Huxley~[9] considered the sum~$S_2$ over a short interval.
The matrices $M$ are conjugates by a fixed matrix of matrices with
small entries. Huxley and Watt [19,20,21] considered the sums
$S_2$ and~$S_3$ with congruence conditions; the matrices $M$ lie
in a congruence subgroup of the modular group.

These improvements relate to the largest range of sums in the
approximate functional equation for $\zeta (s)$, or for the
Fourier series in the lattice point problems. The methods become
weaker for smaller ranges, when the parameters have different
sizes. P.~Sargos [40] gave a variation of the method which works well
at $n = t^\alpha $ in the approximate functional equation for
$\zeta (s)$, with $\alpha = 0.4$, and Huxley and Kolesnik
gave two versions [16,17] of an iterative method which
works well near $\alpha  = 0.42$.

The current best results by this method are
$$
\eqalign{ \zt   \;&\ll_\e\;  t^{32/205 +
\varepsilon }\qquad(32/205 = 0.15609\ldots\,),
 \cr
\Delta (x)  \;&\ll_\e\;  x^{131/416 + \varepsilon },
\qquad(131/416 = 0.314903\ldots\,), \cr
E(T) \;&\ll_\e\;  T^{131/416 + \varepsilon }, \cr}
$$
due to Huxley [14,13] and Watt (unpublished), respectively.

\vskip1cm
\bigskip\no
{\it The Jutila uniform bound

\noindent
Is something not easily found.

With different cases

In different places,

\noindent
When transforms go round and around.}

\vskip2cm


\Refs
\bigskip
\item{[1]} F.V. Atkinson, The mean value of the Riemann zeta-function,
Acta Math. {\bf81}(1949), 353-376.

\item{[2]} R. Balasubramanian,
An improvement on a theorem of Titchmarsh on the mean square of $|\zt|$,
Proc. London Math. Soc. (3) {\bf36}(1978), 540-576.

\item{[3]} E. Bombieri and H. Iwaniec,
On the order of $\zt$,
Ann. Scuola  Norm. Sup. Pisa Cl. Sci. {\bf(4) 13}(1986), 449-472.

\item{[4]} E. Bombieri and H. Iwaniec,
Some mean value theorems for exponential sums,
Ann. Scuola Norm.  Sup. Pisa Cl. Sci. (4) {\bf 13}(1986), 473-486.

\item{[5]} S.W. Graham and G. Kolesnik,
Van der Corput's method of exponential sums,
LMS Lecture Note Series {\bf126}, {\it Cambridge University Press},
Cambridge etc., 1991.

\item{[6]} D.R. Heath-Brown,
Mean value theorems for the Riemann zeta-function,
S\'eminaire de Th\'eorie des Nombres, S\'em. Delange-Pisot-Poitou,
Paris 1979/1980, {\it Birkh\"auser Verlag} 1981, pp. 123-134.

\item{[7]} D.R. Heath-Brown and M. N. Huxley,
Exponential sums with a difference,
Proc. London Math. Soc. (3){\bf 61}(1990), 227-250.

\item{[8]} M.N. Huxley,
Exponential sums and lattice points,
Proc. London Math. Soc. (3) {\bf 60}(1990), 471-502,
Corrigenda ibid. (3) 66, 70.

\item{[9]} M.N. Huxley,
A note on short exponential sums,
in ``Proc. Amalfi Conf. Analytic Number Theory" (eds. E.Bombieri et al.),
University of Salerno, Salerno, 1992, pp.
217-229.

\item{[10]} M.N. Huxley,
Exponential sums and the Riemann zeta function IV,
Proc. London Math. Soc. (3) {\bf 66}(1993), 1-40.

\item{[11]} M.N. Huxley,
 The mean lattice point discrepancy,
Proc. Edinburgh Math. Soc. {\bf 38}(1995), 523-531.

\item{[12]}  M.N. Huxley,
Area, Lattice Points and Exponential Sums,
{\it Oxford Science Publications, Clare-\break ndon Press}, Oxford, 1996.

\item{[13]} M.N. Huxley,
Exponential sums and lattice points III,
Proc. London Math. Soc., (3) {\bf 87}(2003), 591-609,

\item{[14]} M.N. Huxley,
Exponential sums and the Riemann zeta function V,
Proc. London Math. Soc. (3) {\bf 90}(2005), 1-41.

\item{[15]} M.N. Huxley and G. Kolesnik,
Exponential sums and the Riemann zeta function III,
Proc. London Math. Soc. (3) {\bf 62}(1991), 449-468,
Corrigenda ibid. (3) {\bf 66}, 302.

\item{[16]} M.N. Huxley and G. Kolesnik,
Exponential sums with a large second derivative,
preprint 1994.

\item{[17]} M.N. Huxley and G. Kolesnik,
Exponential sums with a large second derivative,
in Number Theory in Memory of Kustaa Inkeri, {\it De Gruyter, Berlin}
2001, 131-144.

\item{[18]}
M.N. Huxley and N. Watt,
Exponential sums and the Riemann zeta function,
Proc. London Math. Soc. (3) {\bf 57}(1988), 1-24.

\item{[19]}
M.N. Huxley and N. Watt,
The number of ideals in a quadratic field,
Proc. Indian Acad. Sci. (Math. Sci.) {\bf 104}(1994), 157-165.

\item{[20]}
M.N. Huxley and N. Watt,
Congruence families of exponential sums,
in ``Analytic Number Theory" (ed. Y. Motohashi), {\it Cambridge University Press}
1997, 127-138.

\item{[21]}
M.N. Huxley and N. Watt,
The number of ideals in a quadratic field II,
Israel J. Math. {\bf 120}(2001), 125-153.

\item{[22]}
A. Ivi\'c, The mean values of the Riemann zeta-function, Tata
Institute of Fundamental Research, Lecture Notes {\bf82}, Bombay
1991 (distr. {\it Springer Verlag, Berlin etc.}), 363 pp.

\item{[23]} A. Ivi\'c, The Riemann zeta-function,
{\it Dover} 2nd ed., Mineola (New York), 2003.

\item{[24]} A. Ivi\'c,
The circle and divisor problem,
Bulletin CXXIX de l'Acad\'emie Serbe des Sciences et des
Arts - 2004, Classe des Sciences math\'ematiques et naturelles,
Sciences math\'ematiques No. {\bf29}, pp. 79-83.

\item{[25]} A. Ivi\'c,
On the Riemann zeta-function and the divisor problem,
Central European J. Math. {\bf(2)(4)} (2004), 1-15, and II, ibid.
{\bf(3)(2)} (2005), 203-214.

\item{[26]} A. Ivi\'c,
On the mean square of the zeta-function and the divisor problem,
 Ann. Acad. Scien. Fennicae Math.{\bf32}(2007), 1-9.

\item{[27]}  H. Iwaniec and C.J.Mozzochi,
On the divisor and circle problems,
J. Number Theory {\bf 29}(1989), 60-93.

\item{[28]} M. Jutila,
Riemann's zeta-function and the divisor problem I, II,
Arkiv Math. {\bf21}(1983), 75-96 and ibid. {\bf31}(1993), 61-70.

\item{[29]} M. Jutila, On a formula of Atkinson,
in ``Coll. Math. Sci. J\'anos Bolyai 34", Topics in Classical Number
Theory, Budapest 1981, {\it North-Holland}, Amsterdam, 1984, pp.
807-823.

\item{[30]} M. Jutila,
Lectures on a Method in the Theory of Exponential Sums,
Tata Institute Lectures in Maths. and Physics 80 (Springer, Bombay) 1987.

\item{[31]} M. Jutila, Mean value estimates for exponential sums,
in ``Number Theory, Ulm 1987", LNM {\bf1380}, Springer Verlag,
Berlin etc. 1989, 120-136.
and II, in Archiv Math. {\bf55}(1990), 267-274.

\item{[32]} M. Jutila, Transformations of exponential sums,
in ``Proc. Amalfi Conf. Analytic Number Theory" (eds. E.Bombieri et al.),
University of Salerno, Salerno, 1992, pp. 263-270.

\item{[33]}
M. Jutila,
The spectral mean square of Hecke $L$-functions on the critical line,
Publ. Inst. Math., Nouv. S\'er. vol. {\bf76(90)} (2004), 41-55.

\item{[34]} A. Kaczorowski and A. Perelli,
The Selberg class: a survey,
in ``Number Theory in Progress, Proc. Conf. in honour of A. Schinzel
(K. Gy\"ory et al. eds)", {\it de Gruyter}, Berlin, 1999, pp. 953-992.

\item{[35]} Y.-K. Lau and K.-M. Tsang,
Omega result for the mean square of the Riemann zeta-function,
Manuscripta Math. {\bf117}(2005), 373-381.

\item{[36]}
J. E. Littlewood,
Quelques cons\'equences de l'hypoth\`ese que la fonction $\zeta (s)$
de Riemann n'a pas des z\'eros dans le demi-plan $R(s) > 1/2$,
Comptes Rendues Acad. Sci. Paris 154 (1912), 263-266.

\item{[37]} T. Meurman,
A generalization of Atkinson's formula to $L$-functions,
Acta Arith. {\bf47}(1986), 351-370.

\item{[38]} H.L. Montgomery,
Topics in Multiplicative Number Theory,
LNM 227, {\it Springer Verlag}, Berlin etc., 1971.

\item{[39]}
W.G. Nowak,
On the order of the lattice rest of a convex planar domain,
Math. Proc. Cam. Philos. Soc. {\bf 98}(1985), 1-4.

\item{[40]}
P. Sargos,
Points entiers au voisinage d'une courbe, sommes trigonom\'etriques
courtes et paires d'exposants,
Proc. London Math. Soc. (3){\bf 70}(1995), 285-312.

\item{[41]}
E.C. Titchmarsh,
The Theory of the Riemann Zeta-function,
{\it Oxford Science Publications} 2nd ed., Oxford, 1986.

\item{[42]} K.-M. Tsang,
Mean square of the remainder term in the Dirichlet divisor problem II,
Acta Arith. {\bf71}(1995), 279-299.

\item{[43]} N. Watt,
A problem on semicubical powers,
Acta Arith. {\bf 52}(1989), 119-140.

\item{[44]} N. Watt,
An elementary treatment of a general Diophantine problem,
Ann. Scuola Norm. Sup. Pisa Cl. Sci. (4) {\bf 15}(1988), 603-614.

\item{[45]} N. Watt,
Exponential sums and the Riemann zeta function II,
J. London Math. Soc. {\bf 39}(1989), 385-404.

\item{[46]} N. Watt,
On differences of semicubical powers,
Monats. Math. {\bf 141}(2004), 45-81.

\item{[47]} J. R. Wilton,
Vorono\"{\i}'s summation formula,
Quart. J. Maths. Oxford (3) (1932), 26-32.

\vskip1cm
\endRefs

\enddocument

\bye